\documentclass[10pt]{amsart}

\usepackage{amsthm} 
\usepackage{amsmath} 
\usepackage{amssymb} 
\usepackage{psfrag}
\usepackage{graphics,pst-plot,graphicx,verbatim}

\DeclareMathOperator{\Mod}{Mod}
\DeclareMathOperator{\Out}{Out}
\DeclareMathOperator{\Aut}{Aut}
\DeclareMathOperator{\Sp}{Sp}
\DeclareMathOperator{\SL}{SL}

\newcommand{\Z}{\mathbb{Z}}

\begin{document}

\title{Dehn twists have roots}

\author{Dan Margalit}
\author{Saul Schleimer}

\address{\hskip-\parindent
  Dan Margalit \\ Department of Mathematics \\
  503 Boston Ave \\
  Tufts University \\
  Medford, MA 02155}
\email{dan.margalit@tufts.edu}

\address{\hskip-\parindent
  Saul Schleimer \\ Department of Mathematics\\
  University of Warwick\\
  Coventry, CV4 7AL, UK}
\email{s.schleimer@warwick.ac.uk}

\thanks{This work is in the public domain.}

\date{\today}


\maketitle

Let $S_g$ denote a closed, connected, orientable surface of genus
$g$, and let $\Mod(S_g)$ denote its mapping class group, that is, the
group of homotopy classes of orientation preserving homeomorphisms of
$S_g$.

{\bf Fact.} If $g \geq 2$, then every Dehn twist in $\Mod(S_g)$ has a
nontrivial root.

It follows from the classification of elements in $\Mod(S_1) \cong \SL(2,\Z)$ that Dehn twists are primitive in the mapping class group of the torus.

For Dehn twists about separating curves, the fact is well-known:
if $c$ is a separating curve then a square root of the Dehn
twist $T_c$ is obtained by rotating the subsurface of $S_g$ on one
side of $c$ through an angle of $\pi$.  In the case of nonseparating
curves, the issue is more subtle.  We give two (equivalent)
constructions of roots below.

\paragraph{\bf Geometric construction.}
Fix $g \geq 2$.  Let $P$ be a regular $(4g-2)$-gon. Glue opposite
sides to obtain a surface $T \cong S_{g-1}$.  The rotation of $P$
about its center through angle $2\pi g/(2g-1)$
induces a periodic map $f$ of $T$.  Notice that $f$ fixes the points $x,y \in T$ that are the images of the vertices of $P$.  Let $T'$ be the surface obtained from $T$ by removing small open disks centered at $x$ and $y$.  Define $f' = f|T'$.

Let $A$ and $B$ be annular neighborhoods of the boundary components of
$T'$.  Modify $f'$ by an isotopy supported in $A \cup B$ so that
\begin{list}{$\cdot$}{
\setlength{\leftmargin}{.5in}
\setlength{\rightmargin}{.5in}
\setlength{\parsep}{0.5ex plus .2ex minus 0ex}
\setlength{\itemsep}{0.2ex plus 0.2ex minus 0ex}
}
\item $f'|\partial T'$ is the identity,
\item $f'|A$ is a $g/(2g-1)$--left Dehn twist, and
\item $f'|B$ is a $(g-1)/(2g-1)$--right Dehn twist.
\end{list}

Identify the two components of $\partial T'$ to obtain a surface $S
\cong S_g$ and let $h : S \to S$ be the induced map.  Then $h^{2g-1}$
is a left Dehn twist along the gluing curve, which is nonseparating.

\paragraph{\bf Algebraic construction.}
Let $c_1, \dots, c_k$ be curves in $S_g$ where $c_i$ intersects $c_{i+1}$ once for each $i$, and all
other pairs of curves are disjoint.  If $k$ is odd, then a regular neighborhood of $\cup
c_i$ has two boundary components, say, $d_1$ and $d_2$, and we have a relation in $\Mod(S_g)$ as follows:
\[ (T_{c_1}^2 T_{c_2} \cdots T_{c_k})^k = T_{d_1} T_{d_2}. \]
This relation comes from the Artin group of type $B_n$, in particular, the factorization of the central element in terms of standard generators \cite{bs}.   
In the case $k=2g-1$, the curves $d_1$ and $d_2$ are isotopic nonseparating curves; call
this isotopy class $d$.  Using the fact that $T_d$ commutes with each $T_{c_i}$, we see that
\[ [(T_{c_1}^2 T_{c_2} \cdots T_{c_{2g-1}})^{1-g} T_d]^{2g-1} = T_d. \]

\paragraph{\bf Other roots.}
All roots of Dehn twists are obtained in a similar way.  That is, if
$f$ is a root of a Dehn twist $T_d$ then the canonical reduction
system for $f$ is $d$ \cite{blm}. By the Nielsen--Thurston classification
for surface homeomorphisms \cite{wpt}, if we cut the surface along $d$, then $f$
restricts to a finite order element.

\paragraph{\bf Roots of half-twists.}  Let $S_{0,2g+2}$ be the sphere with $2g+2$ punctures (or cone points) and let $d$ be a curve in $S_{0,2g+2}$ with 2 punctures on one side and $2g$ on the other.  On the side of $d$ with 2 punctures, we perform a left half-twist, and on the other side we perform a $(g-1)/(2g-1)$--right Dehn twist by arranging the punctures so that one puncture is in the middle, and the other punctures rotate around this central puncture.  The $(2g-1)^{\textrm{\tiny st}}$ power of the composition is a left half-twist about $d$.  Thus, we have roots of half-twists in $\Mod(S_{0,2g+2})$ for $g \geq 2$.  There is a 2-fold orbifold covering $S_g \to S_{0,2g+2}$ where the relation from our algebraic construction above descends to this relation in $\Mod(S_{0,2g+2})$.  A slight generalization of this construction gives roots of half-twists in any $\Mod(S_{0,n})$ with $n \geq 5$.

\paragraph{\bf Roots of elementary matrices.}
If we consider the map $\Mod(S_g) \to \Sp(2g,\Z)$ given by the action of $\Mod(S_g)$ on $H_1(S_g,\Z)$, we
also see that elementary matrices in $\Sp(2g,\Z)$ have roots; for instance, we have
\[
\left(
\begin{array}{rrrr}
1&0&0&1\\
0&1&0&0\\
0&1&-1&1\\
0&1&-1&0\\
\end{array}\right)^3
=
\left(
\begin{array}{rrrr}
1&1&0&0\\
0&1&0&0\\
0&0&1&0\\
0&0&0&1\\
\end{array}\right)_.
\]
By stabilizing, we obtain cube roots of elementary matrices in
$\Sp(2g,\Z)$ for $g \geq 2$.

\paragraph{\bf Roots of Nielsen transformations.}
Let $F_n$ denote the free group generated by $x_1, \dots, x_n$, let $\Aut(F_n)$ denote the group of automorphisms of $F_n$, and
assume $n \geq 2$.  A Nielsen transformation in $\Aut(F_n)$ is an element conjugate to the one given by
$x_1\mapsto x_1x_2$ and $x_k \mapsto x_k$ for $2 \leq k \leq n$.  The
following automorphism is the square root of a Nielsen transformation in
$\Aut(F_n)$ for $n \geq 3$.
\[
\begin{array}{rcl}
x_1 & \mapsto & x_1x_3 \\ x_2 & \mapsto & x_3^{-1}x_2x_3 \\ x_3 &
\mapsto & x_3^{-1}x_2 \\
\end{array}
\]
Taking quotients, this gives a square root of a Nielsen transformation in $\Out(F_n)$ and, multiplying by $-\textrm{Id}$,  a square root of an elementary matrix in $\SL(n, \Z)$, $n \geq 3$.  Finally, our roots of Dehn twists in $\Mod(S)$ can be modified to work
for punctured surfaces, thus giving ``geometric'' roots of Nielsen transformations
in $\Out(F_n)$.

\bibliographystyle{plain}
\bibliography{roots}

\end{document}